\input amssym
\mag=\magstep1

\def\Stab{\operatorname{Stab}}
\def\Hom{\operatorname{Hom}}

\def\det{\operatorname{det}}
\def\ord{\operatorname{ord}}
\def\Ind{\operatorname{Ind}}
\def\Aut{\operatorname{Aut}}
\def\Gal{\operatorname{Gal}}
\def\Kum{\operatorname{Kum}}
\def\Supp{\operatorname{Supp}}
\documentstyle{amsppt}
\document
\topmatter
\author
Tomohide Terasoma
\endauthor
\title
Monodromy weight filtration is independent of $l$
\endtitle
\endtopmatter

\heading
\S 1 Weight monodromy filtration for nilpotent monodromy
\endheading

  In this section, we recall several fundamental properties of
monodromy weight filtration for nilpotent monodromy and prove some
approximation theory.
 Let $K_0$ a local field of characteristic $p$ with a finite residue
field $\kappa = \bold F_q$.  The integer ring of $K_0$ is denoted by $R_0$.
Let $\bold F$ be an algebraic closure
of $\bold F_q$ and $R = R_0 \otimes_{\bold F_q}\bold F$.
The fraction field of $R$ is denoted by $K$ and we use the notations
$$
\align
& Spec (R_0) = S_0, Spec(K_0) = \eta_0, Spec(\kappa_0) = s_0, \\
& Spec (R) = S, Spec(K) = \eta, Spec(\kappa) = s. \\
\endalign
$$
The algebraic closure of $K$ is denoted by $\bar K$ and $Spec(\bar K)$ is
denote by $\bar\eta$.
Let $\Cal X_0$ be a scheme over $S_0$ whose genric geometric fiber
$X_{\bar\eta} = \Cal X_0 \times_{S_0} \bar\eta$ is smooth.
Denote by $\Cal X$ 
the fiber product $\Cal X_0 \times_{\bold F_q}\bold F$.  Then by 
a theorem of Grothendieck,
the action of the Galois group $I = \Gal (\bar K / K)$ on
$H^i(X_{\bar \eta}, \bold Q_l)$ is quasi unipotent.
Let $J$ be an open compact subgroup of $I$ whose action on
$H^i(X_{\bar \eta}, \bold Q_l)$ is unipotent. The
action of $J$ factors through the maximal tame quotient $J^t$ of $J$.
Let $U$ a topological generator of $J^t \simeq \hat\bold Z'(1)$
and $N$ be the logarithm of the action of $U$ on 
$H^i(X_{\bar\eta}, \bold Q_l)$.
In [D], Deligne introduced an increasing filtration $W{\bullet}$ 
with the following properties:
\roster
\item
$N(W_kH^i(X_{\bar\eta}, \bold Q_l)) 
\subset W_{k-2}H^i(X_{\bar\eta}, \bold Q_l)$
\item
The induced map $N : Gr_k^W(H^i(X_{\bar\eta}, \bold Q_l)) \to 
Gr_{k-2}^W(H^i(X_{\bar\eta}, \bold Q_l))$ induces an isomorphism;
$$
N^k : Gr_k^W(H^i(X_{\bar\eta}, \bold Q_l)) \overset{\simeq}\to\to 
Gr_{-k}^W(H^i(X_{\bar\eta}, \bold Q_l)).
$$
\endroster
We define the primitive part $P_k(H^i(X_{\bar\eta}, \bold Q_l))$ of 
$Gr_k^W(H^i(X_{\bar\eta}, \bold Q_l))$ 
($k \geq 0$) by the kernel of $N^{k+1}$.
Let $G$ be an open comact subgroup of $\Gal (\bar K/K_0)$
such that $J = G \cap \Gal (\bar K_0/K)$ acts
$H^i(X_{\bar\eta}, \bold Q_l)$ nilpotently.  Let $L_0$ be the
corresponding extension of $K_0$ and $\bold F_{q'}$
be the residue field.
Then we have the exact sequence
$$
1 \to J \to G \to \Gal (\bold F/ \bold F_{q'})\to 1.
$$
Since the action of $J$ on
the the associated graded vector space $Gr_k^W(H^i(X_{\bar\eta}, \bold Q_l))$
is trivial, $\Gal (\bold F/\bold F_{q'})$ acts on this vector space.

Before studying the action of $\Gal (\bold F/\bold F_{q'})$ on
$Gr_k^W(H^i(X_{\bar\eta}, \bold Q_l))$,
we try to approximate the local situation by the global situation.
Let $g_0 : Y_0 \to C_0$ be a projective generic geometrically smooth
morphism of relative dimension $n$, where
$C_0$ is a projective curve over $\bold F_{q}$. Let
$p$ be an $\bold F_{q}$ valued point of $C_0$. 

We denote by $O_p$ and $K_p$ by the completion of the strucutre 
sheaf of $C_0$ at the point $p$ and its quotient field respectively.
Let $\bar K_p$ be the algebraic clusure of $K_p$,
and $\bar\xi = Spec (\bar K_p)$.
We denote by $Y_{\bar\xi}$ the base extension of $Y_0$ to $\xi$.

\proclaim{Lemma 1.1}
Let $f_0 :\Cal X_0 \to S_0$ be as above
and $l_1, \dots , l_s$ be finite set of primes different from $p$.  
There exist 
a projective smooth curve $C_0$ over $\bold F_q$, and
a projective geometric generically smooth morphism $Y_0 \to C_0$ 
of relative dimension $n$
with the following property:
\roster
\item
There exists an isomorphism 
$$H^i(Y_{\bar\xi}, \bold Q_{l_i}) \simeq
H^i(X_{\bar\eta}, \bold Q_{l_i})
\tag{1.1}
$$
for all $i = 1, \dots, s$.
\item 
There exists an the isomorhpism 
$$
\CD 
\Gal(\bar K_0/K_0)/P_{K_0,i} @>>> \Gal(\bar\bold F_{q}/\bold F_{q}) \\
@V{\overset{\iota}\to\simeq}VV @VV{=}V \\
\Gal (\bar K_p/K_p)/P_{K_p,i} @>>> 
\Gal (\bar\bold F_{q}/\bold F_{q}), \\
\endCD
$$
which is compatible with the isomorphism of (1.1),
where $P_{K_0,i}$ and $P_{K_p,i}$ is the kernel
of 
$$
\Gal(\bar K_0/K) \to\Aut (H^i(X_{\bar\eta}, \bold Q_{l_i}))
$$
and
$$
\Gal(\bar K_p/K_p\otimes_{\bold F_q}\bold F) 
\to\Aut (H^i(Y_{\bar\xi}, \bold Q_{l_i}))
$$
respectively. 
\endroster
\endproclaim
\demo{Proof}
Let us write $R_0= \bold F_q [[t]]$.
By fixing a projective model of $\Cal X_0$ on $S_0$, we can take
a finitely generated sub ring $\bold F_q[t, x_1, \dots, x_m]$
of $\bold F_q[[t]]$ which contains all the coefficients
of $\Cal X_0$.  The element corresponding to $x_i$ is denoted by
$\xi_i(t) \in \bold F_q[[t]]$.  Let 
$M_0 = Spec(\bold F_q[t, x_1, \dots, x_m])$.  By the definition of 
$M_0$ there exists a projective variety $\Cal X_{M_0} \to M_0$ 
of relative dimension $n$ such that the base change of $\Cal X_{M_0}$
by the map $S_0 \to M_0$
is $\Cal X_0$. Let $M = M_0 \otimes_{\bold F_q}\bold F$,
$\Cal X_M = \Cal X_{M_0} \otimes_{\bold F_q}\bold F$, and 
$f_M : \Cal X_M \to M$.  Since it is generically smooth, there exists
an open set $U_0$ of $M_0$ such that
$R^i (f_M)_*$ is a smooth sheaf on $U = U_0\otimes_{\bold F_q}\bold F$.
The complement of $U_0$ in $M_0$ is denoted by $D_0$.
By taking a finite covering of $M_0$ and alterations [J], Theorem 7.3,
we have a proper dominant morphism $\phi :N_0 \to M_0$ 
and an action of a group $G$ on $N_0$ over $M_0$
with the following property:
\roster
\item
Let $Q(N_0)$ and $Q(M_0)$ be the rational function field of $N_0$ and
$M_0$.
The field extension $Q(N_0)^G$ of $Q(M_0)$ is purely inseparable.
\item
The variety $N_0$ is smooth 
and $E_0 =\phi^{-1}(D_0)$ is a $G$-strict normal crossing divisor.
\item
Let $\bold F_{q'}$ be the constant field of $N_0$, 
$N= N_0\otimes_{\bold F_{q'}} \bold F$, 
$E= E_0\otimes_{\bold F_{q'}} \bold F$, 
$f_N : \Cal X_N \to N$ be the base change of $f_M$ by the morphism
$N \to M$.  
Let $\bar\delta$ be a geometric point of $N-E$.
Then the image of the geometric monodromy action of 
$\pi_1(N-E, \bar\delta )$ on the geometric fiber of $R^i(f_N)_*\bold Q_{l_i}$
at $\bar\delta$ is a pro-$l_i$ group.
\endroster
Since the field $Q(N_0)$ is a finite extension of 
$Q(M_0)$, by choosing an embedding of $Q(N_0)$ into the algebraic
clusure of $K_0$, the composite field $L_0 = K_0\cdot Q(N_0)$ is a
finite exitension of $K_0$.  Let $T_0$ be the trait corresponding
to the extension $L_0$. Then, by using valuative criterion of properness,
we have the following diagram.
$$
\CD
T_0 @>>> N_0 \\
@VVV @VVV \\
S_0 @>>> M_0 \\
\endCD
\tag{1.2}
$$
Let $t_0$ be
the image of the closed point of $T_0$ in $N_0$.
Let $G_0$ be the subgroup of $G$ which preserves the image of $T_0$.
By blowing up $t_0$, we get the similar diagram as (1.2).
Repeating the process, we may assume that the image of the trait
$T_0$ intersects only one component $E_{0,i}$ of $E_0$
transversally.  Let $m$ be the maximal ideal corresponding to 
$t_0$ and $\hat N_{0,m}$ be the completion of $N_0$ at $m$ and
$\hat N_{m} = \hat N_{0,m} \times_{\bold F_{q'}} \bold F$.
Since the image of $T_0$ meets $E_{0,i}$ transversally,
the natural homomorphism
$$
\pi_1(T-\{t\})^{(l_i)} \simeq \bold Z_{l_i}(1)
\to \pi_1(\hat N_m - E_i)^{(l_i)}
\simeq \bold Z_{l_i}(1)
$$
is an isomorphism, where $T = T_0 \times_{\bold F_{q'}}\bold F$
and $E_i = E_{0,i} \times_{\bold F_{q'}}\bold F$.
Here we denote by $\pi_1(T-\{t\})^{(l)}$ the maximal pro-$l$
quotient of the fundamental group of $T-\{t\}$.
Now we take a curve $D_0$ in $N_0$ passing through $t_0$ meeting $E_{0,i}$
transversally and equivariant under the action of $G_0$.
Let $Z_0 = (f_{N_0})^{-1}(C_0)$ and take a quotient $Y_0 = Z_0/G_0$ and 
$C_0 = D_0/G_0$ of
$Z_0$ and $D_0$ under the action of $G_0$.  
Let $p_0$ be the image of $t_0$ under the quotient map.
Then the induced
morphisms $Y_0 \to C_0$ and $p_0$ satisfies the required properties.
\enddemo

As an application of this approximation lemma,
we have the following Deligne's theorem for varieties on local fields.
\proclaim{Lemma 1.2 (see also [D])}
The action of $\Gal (\bold F/\bold F_{q'})$ on
$Gr_k^W(H^i(X_{\bar\eta}, \bold Q_l))$ is of pure weight $i+k$.
\endproclaim
\demo{Definition}
The filtration $W$ introduced here is called 
the weight monodromy filtration of 
$H^i(X_{\bar\eta}\bold Q_l)$.
\enddemo

\heading
\S 2 Independence of $l$ in the global situation
\endheading

Let $C_0$ be a proper smooth curve over $\bold F_q$,
$f_0: X_0 \to C_0$ be a projective flat relative $n$-dimensional
morphism with generic geometrically smooth fiber.  
Let $X = X_0 \otimes_{\bold F_q} \bold F$,
$C = C_0 \otimes_{\bold F_q}\bold F$, and $K_0$ and $K$ be the function
field of $C_0$ and $C$ respectively.  Let $s_0$ be a $\bold F_q$-valued point
of $C_0$ and $s$ be the spectrum of the 
algebraic closure of the residue field of $s_0$.
For a rational function
$g \in  K$ and a character 
$\chi :\mu_d(\bold F) \to \bar\bold Q_l^\times$, the associate Kummer sheaf
on $C$ is denoted by $\Kum (g, \chi )$.  
The twist of $\chi$ by an element of Galois group 
$\sigma \in \Gal (\bold F_{q}(\mu_d)/\bold F_q)$ is denoted by
$\chi^{\sigma}$.
Let $\Cal F_0$ be an
etale sheaf on $C_0$ with a finite geometric monodromy.  The restricition
of $\Cal F_0$ to $C$ is denoted by $\Cal F$.
\proclaim{Lemma 2.1}
There exists an element $g \in K_0$, a positive integer $d$,
and an injective character $\chi :\mu_d \to \bar\bold Q_l^\times$
such that
\roster
\item
the support $\Supp (g)$ of the divisor $(g)$ in $C$ is non-empty and
the finite set $\Supp(f) \cup \{ s_0\}$ contains all the points in $C$
where the morphism $f$ or the sheaf $\Cal F$ are not smooth, 
\item
$f(s_0) =1$, 
\item
the order of the class $g$ in $K^\times/K^{d\times}$ is $d$, and
\item
for all $x \in \Supp (g)$ and $\sigma \in \Gal (\bold F_{q'}/\bold F_q)$,
the semi-simplifaication of the 
local monodromy group on
$R^if_*\bar\bold Q_l \otimes \Kum(g, \chi^{\sigma} ) \otimes \Cal F$
at $x$ has no fixed part.
\endroster
\endproclaim
\demo{Proof}
Fix one closed point $x_1$ in $C_0$ different from $s_0$.
By Rieman-Roch theorem, we can choose a rational function $g$ on $C_0$ 
\roster
\item
whose order at $x_1$
is $1$, 
\item
if $f$ or $\Cal F$ are not smooth at a point $x$ different from $s_0$, then
$\ord_x (g) \neq 0$, and
\item
regular at $s_0$ and $g(s_0)=1$.  
\endroster
Choose a sufficiently 
big $d$ and $\chi$ such that $\Kum (g, \chi)$ has the property (4).
Then we get the required $g$, $d$, and $\chi$.
\enddemo
\proclaim{Remark 2.2}
The condition (1) and (4) of Lemma 2.1 implies
$$
H^0_c(U, R^if_*\bar\bold Q_l\otimes \Kum (g, \chi^{\sigma})\otimes \Cal F )=
H^2_c(U, R^if_*\bar\bold Q_l\otimes \Kum (g, \chi^{\sigma})\otimes \Cal F )=0
$$
for all $\sigma \in \Gal (\bold F_{q'}/ \bold F_q)$.
\endproclaim

Now we introduce a covering $\tilde C_0$, $\tilde X_0$ of $C_0$ and
$X_0$.  Let $\bold F_{q'} = \bold F_q(\mu_d)$ and 
$C_1 = C_0 \otimes_{\bold F_q}\bold F_{q'}$.  By Kummer theory,
the $d$-root of $g$ defines a finte cyclic covering $\tilde C_0$ of
$C_1$.  By the composite $\tilde C_0 \to \bold F_{q'} \to \bold F_q$,
$\tilde C_0$ is considered as a curve on $\bold F_q$.  Note that
it is not always geometrically connected. Let $\pi : \tilde C_0 \to C_0$
be the natural projection. Then we have
$$
G = \Aut (\tilde C_0/C_0) \simeq \mu_d \rtimes \Gal (\bold F_{q'}/ \bold F_q).
$$
The induced representation $\Ind_{\mu_d}^G(\chi )$ is denoted by
$\Ind$ and the group ring $\bar\bold Q_l [G]$ is denoted by $A$.
The Kummer sheaf $K_0$ on $C_0$ is defined by
$$
K_0 = \pi_* \bar\bold Q_l \otimes_{A}\Ind.
$$
Then it is easy to see that 
$$
K_0\mid_{C} \simeq \oplus_{\sigma\in\Gal (\bold F_{q'}/ \bold F_q)}
\Kum (g, \chi^{\sigma} ).
$$
Let 
$\tilde X_0 = X_0 \times_{C_0}\tilde C_0$ and 
$\tilde X = \tilde X_0 \otimes_{\bold F_q} \bold F$.  
The projetion $\tilde X_0 \to C_0$ and $\tilde X \to C$ is denoted by
$\tilde f_0$ and $\tilde f$ respectively.
Then the group $G$
acts on $\tilde X_0$ and $\tilde X$.
Let $W_0 = C_0-\{s_0\}$, $U_0 = W_0 -\Supp (g)$ and 
$W = W_0\otimes_{\bold F_q}\bold F$, $U = U_0 \otimes_{\bold F_q}\bold F$.
The natrual inclusions $U \to W$, 
$W \to C$ and $U \to C$ are denoted by $j_1$, $j_2$ and $j_3$ respectively.
\proclaim{Proposition 2.3}
The action of the Frobenius $Frob_{\bold F_q}$ on
$$
H^1(C, (j_2)_*(j_1)_!(R^i\tilde f_*\bar\bold Q_l\otimes_{A}\Ind)\otimes \Cal F)
\tag{2.1}
$$
is pure of weight $i+1$.
\endproclaim
\demo{Proof}
Let $\tilde X^0= X_0\times_{\bold F_{q'}}\bold F$ 
be the connected component of $\tilde X$, and
$\tilde f^0 : \tilde X^0 \to C$ be the natural projection.
Then the 
$\Gal (\bold F/ \bold F_q)$-module (2.1)
is isomorphic to the induced representation of
the $\Gal (\bold F/ \bold F_{q'})$-module
$$
\align
& H^1(C, (j_2)_*(j_1)_!(R^i\tilde f^0_*\bar\bold Q_l
\otimes_{\bold Q_l[\mu_d]}\bar\bold Q_l(\chi ))\otimes \Cal F) 
\tag{2.2}
\\
& \simeq H^1(C, (j_2)_*(j_1)_!(R^i f_*\bar\bold Q_l
\otimes \Kum (g,\chi)\otimes \Cal F). \\
\endalign
$$
Therefore it is enough to prove that (2.2)
is pure of weight $i+1$ under the action of 
$\Gal (\bold F/ \bold F_{q'})$.  By the condition of Lemma 2.1 (4), we have
$$
 (j_1)_!(R^if_*\bar\bold Q_l
\otimes \Kum (g, \chi )\otimes \Cal F) 
 \simeq
 (j_1)_*(R^if_*\bar\bold Q_l
\otimes \Kum (g, \chi)\otimes \Cal F). 
$$
Therefore we have
$$
 (j_2)_*(j_1)_!(R^if_*\bar\bold Q_l
\otimes \Kum (g, \chi )\otimes \Cal F) 
\simeq (j_3)_*(R^if_*\bar\bold Q_l
\otimes \Kum (g, \chi )\otimes \Cal F). 
$$
Since $f$ is projective smooth on $U$, 
$R^if_*\bar\bold Q_l \otimes \Kum(g, \chi ) \otimes \Cal F$ is punctually
pure of weight $i$.
Therefore by the theorem of purity in [D], Theorem 3.2.3, 
we get the proposition.
\enddemo
Now we fix an identification 
$
\bar\bold Q_l \overset{\iota}\to\simeq \bar\bold Q_{l'}
$
for primes $l$ and $l'$ different from $p$.
Let $D_0 \to C_0$ be a Galois covering with the finite Galois group $G$
and $\tau$ be a $\bar\bold Q_l$-valued representation of $G$.
Let $\Cal F_0(\tau)$ be the $l$-adic sheaf associated to the representation
$\tau$.  
The restriction of $\Cal F_0$ to $C$ is denoted by $\Cal F$.
We compare elements of $\bar\bold Q_l$ and $\bar\bold Q_{l'}$
via the isomorphism $\iota$.
\proclaim{Proposition 2.4}
The characteristic polynomial of $Frob_{\bold F_q}$ on
$$
H^1_c(U, R^i\tilde f_*\bar\bold Q_l\otimes_{A}\Ind\otimes \Cal F(\tau))
$$
is indpendent of $l$.
\endproclaim
\demo{Proof}
By Remark 2.2 and Lefschetz trace formula, 
we have the following equality for zeta function of the sheaf 
$R^{i}(\tilde f_0)_* \bar\bold Q_l\otimes_A Ind\otimes \Cal F_0(\tau)$ on $U_0$.
$$
\align
& \det (1 - tFrob_{\bold F_q} \mid 
H_c^1(U, R^{i}\tilde f_* \bar\bold Q_l \otimes_A \Ind \otimes \Cal F(\tau))) 
\\
= &
\prod_{x \in \mid U_0 \mid}
\det (1- t Frob_{x}\mid 
(R^{i}\tilde f_* \bar\bold Q_l \otimes_A \Ind\otimes 
\Cal F(\tau))_{\bar x})^{-1} \\
= &
\prod_{x \in \mid U_0 \mid}
\det (1- t Frob_{x}\mid 
H^{i}((\tilde f)^{-1}(\bar x), \bar\bold Q_l) \otimes_A \Ind\otimes 
\Cal F(\tau)_{\bar x})^{-1} \\
\endalign
$$
Since $\tilde f$ is projective smooth at $x \in U$, the right hand
side is independent of $l$ by the classical Weil conjecture and relation
between Zeta function and the Frobenius action and
the action of $G$ on the $\bold F$ rational points of
$\tilde f^{-1}(x)$.  Therefore we get the proposition.
\enddemo

\heading
\S 3 Proof of the main theorem
\endheading

Let $K_0$ be a local field of characteristic $p>0$ with a
residue field $\bold F_q$.
Let $X_\eta$ be a projective geometrically smooth variety of dimension $n$
over $K_0$.
The variety $X_\eta$ is called globalizable if there exist
projective varieties $C_0$, $X_0$ over $\bold F_q$ of dimension
one and $n+1$, a morphism $g : X_0 \to C_0$, a point 
$p \in C_0(\bold f_q)$ and an isomorphism of local field between $K_0$
and the quotient field $K_p$ of the completion of the strucutre sheaf of
$C_0$ at $p$ such that
the base change of $X_0$ by the morphism $Spec (K_0) \to Spec (K_p) \to
C_0$ is isomorphic to $X_\eta$.
Let $\tau$ be a representation of
$\pi_1 ( C_0)$ with a finite image and $\bar\bold Q_l (\tau)$
be the representation space of $\tau$.
\proclaim{Proposition 3.1}
The characteristic polynomial of the Frobenius action $Frob_{\kappa (s_0)}$
on the sapce $(H^i(X_{\bar\eta}, \bold Q_{l_1})
\otimes \bar\bold Q_l (\tau))^{\Gal (\bar K_p/ K_p)}$ 
is independent of $l$.
\endproclaim
\demo{Proof}
To show the independence of $l$, it is enough to compare two primes 
$l_1$ and $l_2$.  By using Lemma 1.1, we may assume that $X_\eta$
is globalizable.  We use the same notation for $X_0$, $C_0$, $p_0$
and so on.
First we choose a rational function $g$ on $C_0$, an integer $d$, and
a character $\chi : \mu_d(\bold F) \to \bar\bold Q_l^\times$ satisfying
the properties in Lemma 2.1.  We consider varieties $\tilde C_0$, $\tilde C$,
$\tilde X_0$, $\tilde X$ and so on as in \S 2.
Now we consider the exact sequence of sheaves on $C$:
$$
\align
0 & \to (j_3)_!(R^i\tilde f_* \bar\bold Q_l\otimes_A Ind)\otimes \Cal F (\tau)
 \to(j_2)_*((j_1)_!(R^i\tilde f_* \bar\bold Q_l\otimes_A Ind) 
\otimes\Cal F (\tau))
\\
& \to
i_{s}^*((j_3)_*(R^i\tilde f_* \bar\bold Q_l\otimes_A Ind\otimes \Cal F(\tau)))
\to 0, \\
\endalign
$$
where $i_s : \{s\} \to C$ is the natural inclusion.
Taking the cohomology, we have the long exact sequence;
$$
\align
0 & \to
H^0(Spec (K), H^i(\tilde X_{\bar\eta}, \bar\bold Q_l)\otimes_A Ind 
\otimes \bar\bold Q_l(\tau))
 \to
H^1_c(U, R^i\tilde f_* \bar\bold Q_l\otimes_A Ind \otimes \Cal F(\tau)) \\
& \to
H^1( C, (j_2)_*(j_1)_!(R^i\tilde f_* \bar\bold Q_l\otimes_A Ind
\otimes \Cal F(\tau)))
\to 0. \\
\endalign
$$ 
The weight of
$$
H^0(Spec (K), H^i(\tilde X_{\bar\eta}, \bar\bold Q_l)\otimes_A Ind
\otimes \bar\bold Q_l(\tau))
\tag{3.1}
$$
is less than or equal to $i$, 
and that of 
$H^1( C, (j_2)_*((j_1)_!(R^i\tilde f_* \bar\bold Q_l\otimes_A Ind)
\otimes \Cal F (\tau)))$
is purely $i+1$.
Therefore the characteristic polynomial of
the Frobenius action on (3.1) is independent of $l$.
By the property Lemma 2.1 (2) of the choice 
of $g$, (3.1) is isomorphic to
$$
H^0(Spec(K), H^i(X_{\bar\eta} \bar\bold Q_l)
\otimes \bar\bold Q_l(\tau))^{\oplus r},
$$
where $r = [\bold F_q(\mu_d): \bold F_q]$, and we get the theorem.
\enddemo
By using Poincare duality and local duality, we have the following
corollary.
\proclaim{Corollary 3.2}
The characteristic polynomial of the Frobenius action on
\linebreak
$H^1(Spec (K),H^{i}(X_{\bar\eta},\bold Q_l) \otimes \bar\bold Q_l(\tau))$
is independent of $l$.
\endproclaim

Now we recall some properties of the weight monodromy filtration 
$W_{\bullet}$.  Let $G_0$ and $G$ be the absolute Galois groups of 
$K_0$ and $K$ respectively.  There exists an open compact subgroup
$\tilde M^0$ which acts on $H^i( X_{\bar\eta}, \bold Q_l)$ unipotently
and the corresponding quotient of $\tilde M^0$ is denoted by $M^0$. It is 
isomorphic to $\bold Z_l(1)$.  Then the image $M$ of $G$ in 
$\Aut H^i(X_{\bar\eta}, \bar\bold Q_l)$ contains $M^0$ and we may assume
that $M^0$ is a normal subgroup of $M$ by changing $M^0$ sufficiently
small.  Then the quotient group $M/M^0$ acts on $\bold Z_l(1)$
by the conjugation and the corresponding character is denoted by $\beta$.
Then the relation of 
$N : Gr_{k}^W \to Gr_{k-2}^W$ and $g \in M/M^0$ is given by
$Ng = \beta (g)^{-1}gN$.  
Therefore,
by the universal property of the filtration $W$, the action of the 
group $M/M^0$ preserves the filtration $W$ and as a consequence $M/M^0$
acts on the associate graded module 
$Gr_k^W = Gr_k^W(H^i(X_{\bar\eta}, \bar\bold Q_l))$.
Again by the relation of $N$ and $g$,
the group $M/M^0$ acts
on the primitive part $P_k$. The corresponding representation of
$P_k$ is denoted by $\alpha_k$.  
Then as representations of $M/M^0$,
$Gr_k^W$, $H^i(X_{\bar\eta}, \bold Q_l)^{M^0}$,
and the coinvariant $H^i(X_{\bar\eta}, \bold Q_l)_{M^0}$ under the action
of $M^0$ are isomorphic to
$$
\align
Gr_k^W &= \oplus_{m\geq 0}\alpha_{k+2m}\otimes \beta^{m}, \\
H^i(X_{\bar\eta}, \bold Q_l)^{M^0} &= \oplus_{m\geq 0}
\alpha_m\otimes \alpha^m, \\
H^i(X_{\bar\eta}, \bold Q_l)_{M^0} &=\oplus_{m\geq 0}\beta_m. \\
\endalign
$$
It is easy to see that all the intersection
of the kernel of $\alpha$ and $\alpha_k$ corresponds to
the maximal nilpotent subgroup for the action of $M$ on 
$\Aut (H^1(X_{\bar\eta}, \bar\bold Q_l))$.  Changing notation, 
this group is denoted by $M^0$.  
Let $\tilde M^0$ be the inverse image of $M^0$ under the natural map
$G \to M$.
Let $N_0$ be an open compact subgroup of $G_0$ where
the intersection $N_0 \cap G$ is contained in $\tilde M_0$.
The corresponding extension of $K_0$ is denoted by
$L_0$.  The residue field of $L_0$ is denoted by $\bold F_{q'}$.
\proclaim{Theorem 3.3}
\roster
\item
The representation $\alpha_k$ and $\beta$ of $G$ is independent of $l$.
Especially, $\tilde M^0$ does not depend on $l$ and the dimension
of $P_k$ is independent of $l$.
\item
The characteristic polynomial of the Forbenius action $Frob_{\bold F_{q'}}$
on $P_k$ is independent of $l$.
\endroster
\endproclaim
\demo{Proof of (1)}
Let $\tau$ be a finite dimensional irreducible representation witha finte
image $H$.  Then there exists a Galois covering $\pi : D \to C$
and a lifting $q$ of $p$ such that the quotient field $K_{D,q}$ of
the completion of the structure sheaf $\Cal O_D$ of $D$ at $q$ corresponds
to the quotient $H$ of $G$.  Then the curve $C$, $D$ and all the morphism
in $\Gal (D/C)$ is defined over a finite extension $\bold F_{q'}$
of $\bold F_q$.  The model of $C$ and $D$ defined on $\bold F_{q'}$
is denoted by $C_1$ and $C_2$ respectively  Let $C_1'$ be the
covering between $D \to C$ corresponding to the stabilizer $\Stab_q(G)$
of $q$.  Since $H \simeq \Stab_q(G) \simeq \Gal (D_1/C_1')$,
we can consider a sheaf $\Cal F_1'(\tau)$ on $C_1'$.  We apply 
Proposition 3.1 to $\Cal F_1(\tau)$ and $C_1'$.  Then the characteristic 
polynomial of the Frobenius action on
$(\bar\bold Q_l(\tau) \otimes H^i(X_{\bar\eta}, \bold Q_l))^G$ is 
independent of $l$.  Let $\tilde M^0_l$ be the maximal unipotent
subgroup for the representation $H^i(X_{\bar\eta}, \bar\bold Q_l)$.
Then there exists a finite normal subgroup $\tilde M^1_l$
of $\tilde M^0_l$ such that the restriction of $\tau$ to $\tilde M^1_l$
is trivial.  Then the characteristic polynomila of Frobenius action on
$$
\align
& (\bar\bold Q_l (\tau) \otimes H^i(X_{\bar\eta}, \bar\bold Q_l))^G  \\
= &
((\bar\bold Q_l (\tau) \otimes H^i(X_{\bar\eta}, \bar\bold Q_l))^{\tilde M^1_l}
)^{G/\tilde M^1_l} \\
= &
(\bar\bold Q_l (\tau) \otimes (H^i(X_{\bar\eta}, \bar\bold Q_l))^{\tilde M^1_l}
)^{G/\tilde M^1_l} \\
= &
(\bar\bold Q_l (\tau) \otimes (H^i(X_{\bar\eta}, \bar\bold Q_l))^{\tilde M^0_l}
)^{G} \\
= &
\Hom_G (\bar\bold Q_l (\tau^*), 
H^i(X_{\bar\eta}, \bar\bold Q)^{\tilde M^0_l}) \\
\endalign
$$
is independent of $l$.  The dimension of weight $k$-part of the above
equality is nothing but the multiplicity of $\tau^*$ in 
$\alpha_k\otimes\beta^k$.  Since $\tau$ is the arbitrary irreducible 
representation of $G$ with finite image, $\alpha_k \otimes\beta^k$
is independent of $l$.
Similarly, the characteristic polynomial of Frobenius action on
$$
\align
& H^1(Spec(K), H^i(X_{\bar\eta}, \bar\bold Q_l)\otimes \bar\bold Q_l(\tau ))
\\
& =
\Hom_G (\bar\bold Q_l (\tau^*), 
H^i(X_{\bar\eta}, \bar\bold Q_l)_{\tilde M^0_l}) \\
\endalign
$$
is independent of $l$ and we have the $l$-independence of $\alpha_k$
($k \geq 0$).  This proves (1) of Theorem 3.3.  Therefore by the
characterization of $\tilde M^0_l$ gives as above, it is independent of $l$.

Let $D$ be a covering of $C$ and $q \in D$ a lifting
of $p$ such that the field $K_{D,q}$ defined as above corresponds
to the sub group $\tilde M^0$ of $G$.  Then $D$, $C$ are defined over some 
finite extension $\bold F_{q'}$ of $\bold F_q$ and the model 
over $\bold F{q'}$ are denoted by $D_1$ and $C_1$ respectively.  
We may assume $q \in D_1(\bold F_{q'})$.  We apply Proposition 3.1
to $q \in D_1$ and the trivial representation of $\pi_1(D_1)$.
Then the characteristeic polynomial of Frobenius action on 
$H^i(X_{\bar\eta}, \bar\bold Q_l)^{\tilde M^0}$ is independent of $l$.
Considering weight of $P_k$, the caracteristic polynomial of Frobenius
on $P_k$ is independent of $l$.
\enddemo
\proclaim{Remark 3.4}
\roster
\item
By the independence of $l$ for $P_k$ implies the independence for
the Swan conductor.
\item
Since $\beta$ is independent of $l$, it is order 2 character by
Chevotalev density theorem.
\endroster
\endproclaim

\enddocument